\documentclass{amsart}

\usepackage{amsfonts,amssymb,verbatim,amsmath,amsthm,latexsym,textcomp,amscd}
\usepackage{latexsym,amsfonts,amssymb,epsfig,verbatim}
\usepackage{amsmath,amsthm,amssymb,latexsym,graphics,textcomp}
\usepackage{graphicx}
\usepackage{color}
\usepackage{url}

\input xy
\xyoption{all}

\begin{document}

\newtheorem{theorem}{Theorem}[section]
\newtheorem{prop}[theorem]{Proposition}
\newtheorem{lemma}[theorem]{Lemma}
\newtheorem{cor}[theorem]{Corollary}
\newtheorem{definition}[theorem]{Definition}
\newtheorem{conj}[theorem]{Conjecture}
\newtheorem{rmk}[theorem]{Remark}
\newtheorem{claim}[theorem]{Claim}
\newtheorem{defth}[theorem]{Definition-Theorem}

\newcommand{\boundary}{\partial}
\newcommand{\C}{{\mathbb C}}
\newcommand{\integers}{{\mathbb Z}}
\newcommand{\natls}{{\mathbb N}}
\newcommand{\ratls}{{\mathbb Q}}
\newcommand{\bbR}{{\mathbb R}}
\newcommand{\proj}{{\mathbb P}}
\newcommand{\lhp}{{\mathbb L}}
\newcommand{\tube}{{\mathbb T}}
\newcommand{\cusp}{{\mathbb P}}
\newcommand\AAA{{\mathcal A}}
\newcommand\BB{{\mathcal B}}
\newcommand\CC{{\mathcal C}}
\newcommand\DD{{\mathcal D}}
\newcommand\EE{{\mathcal E}}
\newcommand\FF{{\mathcal F}}
\newcommand\GG{{\mathcal G}}
\newcommand\HH{{\mathcal H}}
\newcommand\II{{\mathcal I}}
\newcommand\JJ{{\mathcal J}}
\newcommand\KK{{\mathcal K}}
\newcommand\LL{{\mathcal L}}
\newcommand\MM{{\mathcal M}}
\newcommand\NN{{\mathcal N}}
\newcommand\OO{{\mathcal O}}
\newcommand\PP{{\mathcal P}}
\newcommand\QQ{{\mathcal Q}}
\newcommand\RR{{\mathcal R}}
\newcommand\SSS{{\mathcal S}}
\newcommand\TT{{\mathcal T}}
\newcommand\UU{{\mathcal U}}
\newcommand\VV{{\mathcal V}}
\newcommand\WW{{\mathcal W}}
\newcommand\XX{{\mathcal X}}
\newcommand\YY{{\mathcal Y}}
\newcommand\ZZ{{\mathcal Z}}
\newcommand\CH{{\CC\HH}}
\newcommand\PEY{{\PP\EE\YY}}
\newcommand\MF{{\MM\FF}}
\newcommand\RCT{{{\mathcal R}_{CT}}}
\newcommand\PMF{{\PP\kern-2pt\MM\FF}}
\newcommand\FL{{\FF\LL}}
\newcommand\PML{{\PP\kern-2pt\MM\LL}}
\newcommand\GL{{\GG\LL}}
\newcommand\Pol{{\mathcal P}}
\newcommand\half{{\textstyle{\frac12}}}
\newcommand\Half{{\frac12}}
\newcommand\Mod{\operatorname{Mod}}
\newcommand\Area{\operatorname{Area}}
\newcommand\ep{\epsilon}
\newcommand\hhat{\widehat}
\newcommand\Proj{{\mathbf P}}
\newcommand\U{{\mathbf U}}
 \newcommand\Hyp{{\mathbf H}}
\newcommand\D{{\mathbf D}}
\newcommand\Z{{\mathbb Z}}
\newcommand\R{{\mathbb R}}
\newcommand\Q{{\mathbb Q}}
\newcommand\E{{\mathbb E}}
\newcommand\til{\widetilde}
\newcommand\length{\operatorname{length}}
\newcommand\tr{\operatorname{tr}}
\newcommand\gesim{\succ}
\newcommand\lesim{\prec}
\newcommand\simle{\lesim}
\newcommand\simge{\gesim}
\newcommand{\simmult}{\asymp}
\newcommand{\simadd}{\mathrel{\overset{\text{\tiny $+$}}{\sim}}}
\newcommand{\ssm}{\setminus}
\newcommand{\diam}{\operatorname{diam}}
\newcommand{\pair}[1]{\langle #1\rangle}
\newcommand{\T}{{\mathbf T}}
\newcommand{\inj}{\operatorname{inj}}
\newcommand{\pleat}{\operatorname{\mathbf{pleat}}}
\newcommand{\short}{\operatorname{\mathbf{short}}}
\newcommand{\vertices}{\operatorname{vert}}
\newcommand{\collar}{\operatorname{\mathbf{collar}}}
\newcommand{\bcollar}{\operatorname{\overline{\mathbf{collar}}}}
\newcommand{\I}{{\mathbf I}}
\newcommand{\tprec}{\prec_t}
\newcommand{\fprec}{\prec_f}
\newcommand{\bprec}{\prec_b}
\newcommand{\pprec}{\prec_p}
\newcommand{\ppreceq}{\preceq_p}
\newcommand{\sprec}{\prec_s}
\newcommand{\cpreceq}{\preceq_c}
\newcommand{\cprec}{\prec_c}
\newcommand{\topprec}{\prec_{\rm top}}
\newcommand{\Topprec}{\prec_{\rm TOP}}
\newcommand{\fsub}{\mathrel{\scriptstyle\searrow}}
\newcommand{\bsub}{\mathrel{\scriptstyle\swarrow}}
\newcommand{\fsubd}{\mathrel{{\scriptstyle\searrow}\kern-1ex^d\kern0.5ex}}
\newcommand{\bsubd}{\mathrel{{\scriptstyle\swarrow}\kern-1.6ex^d\kern0.8ex}}
\newcommand{\fsubeq}{\mathrel{\raise-.7ex\hbox{$\overset{\searrow}{=}$}}}
\newcommand{\bsubeq}{\mathrel{\raise-.7ex\hbox{$\overset{\swarrow}{=}$}}}
\newcommand{\tw}{\operatorname{tw}}
\newcommand{\base}{\operatorname{base}}
\newcommand{\trans}{\operatorname{trans}}
\newcommand{\rest}{|_}
\newcommand{\bbar}{\overline}
\newcommand{\UML}{\operatorname{\UU\MM\LL}}
\newcommand{\EL}{\mathcal{EL}}
\newcommand{\tsum}{\sideset{}{'}\sum}
\newcommand{\tsh}[1]{\left\{\kern-.9ex\left\{#1\right\}\kern-.9ex\right\}}
\newcommand{\Tsh}[2]{\tsh{#2}_{#1}}
\newcommand{\qeq}{\mathrel{\approx}}
\newcommand{\Qeq}[1]{\mathrel{\approx_{#1}}}
\newcommand{\qle}{\lesssim}
\newcommand{\Qle}[1]{\mathrel{\lesssim_{#1}}}
\newcommand{\simp}{\operatorname{simp}}
\newcommand{\vsucc}{\operatorname{succ}}
\newcommand{\vpred}{\operatorname{pred}}
\newcommand\fhalf[1]{\overrightarrow {#1}}
\newcommand\bhalf[1]{\overleftarrow {#1}}
\newcommand\sleft{_{\text{left}}}
\newcommand\sright{_{\text{right}}}
\newcommand\sbtop{_{\text{top}}}
\newcommand\sbot{_{\text{bot}}}
\newcommand\sll{_{\mathbf l}}
\newcommand\srr{_{\mathbf r}}
\newcommand\geod{\operatorname{\mathbf g}}
\newcommand\mtorus[1]{\boundary U(#1)}
\newcommand\A{\mathbf A}
\newcommand\Aleft[1]{\A\sleft(#1)}
\newcommand\Aright[1]{\A\sright(#1)}
\newcommand\Atop[1]{\A\sbtop(#1)}
\newcommand\Abot[1]{\A\sbot(#1)}
\newcommand\boundvert{{\boundary_{||}}}
\newcommand\storus[1]{U(#1)}
\newcommand\Momega{\omega_M}
\newcommand\nomega{\omega_\nu}
\newcommand\twist{\operatorname{tw}}
\newcommand\modl{M_\nu}
\newcommand\MT{{\mathbb T}}
\newcommand\Teich{{\mathcal T}}
\renewcommand{\Re}{\operatorname{Re}}
\renewcommand{\Im}{\operatorname{Im}}

\title{Moebius rigidity for compact deformations of negatively curved manifolds}

\author{Kingshook Biswas}
\address{Indian Statistical Institute, Kolkata, India. Email: kingshook@isical.ac.in}

\begin{abstract} Let $(X, g_0)$ be a complete, simply connected Riemannian manifold with sectional curvatures $K_{g_0}$
satisfying $-b^2 \leq K_{g_0} \leq -1$ for some $b \geq 1$. Let $g_1$ be a Riemannian metric on $X$ such that $g_1 = g_0$ outside
a compact in $X$, and with sectional curvatures $K_{g_1}$ satisfying $K_{g_1} \leq -1$. The identity map $id : (X, g_0) \to (X, g_1)$
is bi-Lipschitz, and hence induces a homeomorphism between the boundaries at infinity of $(X, g_0)$ and $(X, g_1)$, which we denote by
$\hat{id}_{g_0, g_1} : \partial_{g_0} X \to \partial_{g_1} X$. We show that if the boundary map $\hat{id}_{g_0, g_1}$ is Moebius
(i.e. preserves cross-ratios), then it extends to an isometry $F : (X, g_0) \to (X, g_1)$.
\end{abstract}

\bigskip

\maketitle

\tableofcontents

\section{Introduction}

\medskip

In various rigidity problems for negatively curved spaces, the interplay between the geometry of the space and the geometry of its
boundary at infinity plays a prominent role. For a CAT(-1) space $X$, there is a positive function called the {\it cross-ratio} defined
on the space of quadruples of distinct points in the boundary $\partial X$, and a well-known problem asks whether the cross-ratio in fact
determines the space up to isometry. More precisely, if $f : \partial X \to \partial Y$ is a Moebius homeomorphism between boundaries of
CAT(-1) spaces $X, Y$ (i.e. a homeomorphism which preserves cross-ratios), then the question is whether $f$ extends to an isometry
$F : X \to Y$. It is a classical result that this holds when $X = Y = \mathbb{H}^n$, the real hyperbolic space, a fact which is often
used in rigidity theorems for hyperbolic manifolds, for example in the Mostow Rigidity theorem \cite{mostow-ihes}. More generally, Bourdon \cite{bourdon2}
showed that if $X$ is a rank one symmetric space of noncompact type with the metric normalized such that the maximum of the sectional curvatures
equals $-1$, and $Y$ is any CAT(-1) space, then any Moebius embedding $f : \partial X \to \partial Y$ extends to an isometric embedding
$F : X \to Y$. For general CAT(-1) spaces $X, Y$, the problem remains open.

\medskip

We should remark that one of the main motivations for
studying this problem is its relation to the {\it marked length spectrum rigidity} conjecture of Burns and Katok, which asks whether
two closed negatively curved manifolds $X, Y$ with the same marked length are necessarily isometric. Otal \cite{otal2} and independently Croke 
\cite{croke} proved that marked length spectrum rigidity holds in two dimensions. It is well known that in fact
$X, Y$ have the same marked length spectrum if and only if there is an equivariant Moebius map between the boundaries of the universal covers
$f : \partial \tilde{X} \to \partial \tilde{Y}$, so a positive answer to the problem of extending Moebius maps to isometries would
also give a solution to the marked length spectrum rigidity problem (see \cite{otal1}). Equality of marked length spectra is also known to be
equivalent to existence of a homeomorphism between the unit tangent bundles $\phi : T^1 X \to T^1 Y$
conjugating the geodesic flows of $X, Y$ (\cite{hamenstadt1}). Proofs of these equivalences may be found in \cite{biswas3}, section 5. 
We remark that in related work Beyrer, Fioravanti and Incerti-Medici have constructed a cross-ratio on the Roller boundary of any CAT(0) cube complex, and 
have shown that any cross-ratio preserving bijection between geodesically complete cube complexes admits a unique extension to an 
isomorphism of cube-complexes, and have also proved a version of marked length spectrum rigidity for group actions on CAT(0) cube complexes 
\cite{beyrerandall}.  

\medskip

In \cite{biswas3}, it was shown that a Moebius homeomorphism between the boundaries of proper, geodesically complete CAT(-1) spaces
extends to a $(1, \log 2)$-quasi-isometry between the spaces. For $X, Y$ complete, simply connected manifolds of pinched negative curvature
$-b^2 \leq K \leq -1$, this result was refined in \cite{biswas5} to show that the extension may be taken in this case to be a $(1, (1 - 1/b)\log 2)$-quasi-isometry.
In fact the quasi-isometric extension of \cite{biswas3} and \cite{biswas5} was shown to be given by a certain natural extension of Moebius maps called the
{\it circumcenter extension}, which is natural with respect to composition with isometries.
In \cite{biswas6}, it was shown that if $f : \partial X \to \partial Y$ and $g : \partial Y \to \partial X$ are mutually inverse
Moebius homeomorphisms between boundaries of complete, simply connected manifolds $X, Y$ of pinched negative curvature
$-b^2 \leq K \leq -1$, then the circumcenter extensions $F : X \to Y$ and $G : Y \to X$ of $f, g$ are $\sqrt{b}$-bi-Lipschitz homeomorphisms which are inverses
of each other.

\medskip

In the present article we consider compactly supported deformations of the metric on a complete, simply connected manifold
$(X, g_0)$ of pinched negative
curvature $-b^2 \leq K_{g_0} \leq -1$, i.e. we consider metrics $g_1$ on $X$ such that $g_1 = g_0$ outside a compact in $X$, and such that
the sectional curvature of $g_1$ is bounded above by $-1$. The identity map $id : (X, g_0) \to (X, g_1)$ is clearly bi-Lipschitz, hence it induces a homeomorphism between boundaries which we denote by $\hat{id}_{g_0, g_1} : \partial_{g_0} X \to \partial_{g_1} X$, and the problem in this context becomes the following: if $\hat{id}_{g_0, g_1}$
is Moebius, then does it extend to an isometry $F : (X, g_0) \to (X, g_1)$? Partial results for this problem were obtained in \cite{biswas4}, where local and
infinitesimal versions of the problem were considered, namely metrics $g_1$ such that the $C^2$ norm $||g_0 - g_1||_{C^2}$ is small, and one-parameter
families of metrics $(g_t)_{0 \leq t \leq 1}$, and in both cases it was shown that if the boundary maps are Moebius then they extend to isometries.
Our main theorem below gives a complete solution to this problem:

\medskip

\begin{theorem} \label{mainthm} Let $(X, g_0)$ be a complete, simply connected manifold of pinched negative curvature $-b^2 \leq K_{g_0} \leq -1$.
Let $g_1$ be a metric on $X$ such that $g_1 = g_0$ outside a compact in $X$, and such that the sectional curvature of $g_1$ satisfies $K_{g_1} \leq -1$.
Let $\hat{id}_{g_0, g_1} : \partial_{g_0} X \to \partial_{g_1} X$ denote the homeomorphism between boundaries induced by the identity map
$id : (X, g_0) \to (X, g_1)$. Suppose $\hat{id}_{g_0, g_1}$ is Moebius. Then the circumcenter extension of $\hat{id}_{g_0, g_1}$ is an
isometry $F : (X, g_0) \to (X, g_1)$.
\end{theorem}

\medskip

The key to the proof of the above theorem is a further study of properties of the circumcenter extension. In section 2 we briefly recall some facts about
Moebius maps, geodesic conjugacies and circumcenter extensions. In section 3 we prove the results about the circumcenter extension which are used in the
proof of the main theorem, while section 4 is devoted to the proof of the main theorem.

\medskip

\section{Preliminaries}

\medskip

We give below a brief outline of the background on Moebius maps which we will be needing, for details and proofs of the
assertions below the reader is referred to \cite{biswas3}, \cite{biswas5},
\cite{biswas6}.

\medskip

Let $(Z, \rho_0)$ be a compact metric space of diameter one. The cross-ratio with respect to a metric $\rho$ on $Z$ is the function on quadruples of distinct points in $Z$
defined by
$$
[\xi, \xi', \eta, \eta']_{\rho} := \frac{\rho(\xi, \eta) \rho(\xi', \eta')}{\rho(\xi, \eta') \rho(\xi', \eta)} \ , \xi, \xi', \eta, \eta' \in Z
$$
Two metrics $\rho_1, \rho_2$ on $Z$ are said to be Moebius equivalent if their cross-ratios are equal, $[.,.,.,.]_{\rho_1} = [.,.,.,.]_{\rho_2}$.
A metric $\rho$ on $Z$ is said to be antipodal if it has diameter one and for all $\xi \in Z$ there exists $\eta \in Z$ such that $\rho(\xi, \eta) = 1$.
Assume that the metric $\rho_0$ is antipodal. We then define $\mathcal{M}(Z, \rho_0)$ to be the set of all antipodal metrics $\rho$ on $Z$ which
are Moebius equivalent to $\rho_0$. Then for any two metrics $\rho_1, \rho_2 \in \mathcal{M}(Z, \rho_0)$, there is a positive continuous function on $Z$
called the derivative of the metric $\rho_2$ with respect to the metric $\rho_1$, denoted by $\frac{d\rho_2}{d\rho_1}$, such that
$$
\rho_2(\xi, \eta)^2 = \frac{d\rho_2}{d\rho_1}(\xi) \frac{d\rho_2}{d\rho_1}(\eta) \rho_1(\xi, \eta)^2
$$
for all $\xi, \eta \in Z$. If $\xi$ is not an isolated point of $Z$, then 
$$
\frac{d\rho_2}{d\rho_1}(\xi) = \lim_{\eta \to \xi} \frac{\rho_2(\xi, \eta)}{\rho_1(\xi, \eta)} 
$$
Moreover
$$
\left(\max_{\xi \in Z} \frac{d\rho_2}{d\rho_1}(\xi)\right) \left(\min_{\xi \in Z} \frac{d\rho_2}{d\rho_1}(\xi)\right) = 1
$$
This allows us to define a metric on the set $\mathcal{M}(Z, \rho_0)$ by
$$
d_{\mathcal{M}}(\rho_1, \rho_2) := \max_{\xi \in Z} \log \frac{d\rho_2}{d\rho_1}(\xi)
$$
The metric space $(\mathcal{M}(Z, \rho_0), d_{\mathcal{M}})$ is proper and complete. The following lemma follows from the proof of Lemma 2.6 of \cite{biswas3},
we include a proof for convenience:

\medskip

\begin{lemma} \label{maxminantipodal} For $\rho_1, \rho_2 \in \mathcal{M}(Z, \rho_0)$, let $\xi, \eta \in Z$ be points where $\frac{d\rho_2}{d\rho_1}$ attains its
maximum and minimum values respectively. If $\xi' \in Z$ is such that $\rho_1(\xi, \xi') = 1$, then $\frac{d\rho_2}{d\rho_1}$ attains its minimum at $\xi'$, and
$\rho_2(\xi, \xi') = 1$. If $\eta' \in Z$ is such that $\rho_2(\eta, \eta') = 1$, then $\frac{d\rho_2}{d\rho_1}$ attains its maximum at $\eta'$, and
$\rho_1(\eta, \eta') = 1$.
\end{lemma}

\medskip

\noindent{\bf Proof:} Let $\lambda, \mu > 0$ be the maximum and minimum values of $\frac{d\rho_2}{d\rho_1}$ respectively, then we know that $\lambda \cdot \mu = 1$.
For $\xi' \in Z$ such that $\rho_1(\xi, \xi') = 1$, we have
$$
1 \geq \rho_2(\xi, \xi')^2 = \frac{d\rho_2}{d\rho_1}(\xi) \frac{d\rho_2}{d\rho_1}(\xi') \rho_1(\xi, \xi')^2 \geq \lambda \cdot \mu \cdot 1 = 1
$$
so equality holds in the inequalities above, hence $\frac{d\rho_2}{d\rho_1}(\xi') = \mu$ and $\rho_2(\xi, \xi') = 1$.

\medskip

For $\eta' \in Z$ such that $\rho_2(\eta, \eta') = 1$, we have
$$
1 = \rho_2(\eta, \eta')^2 = \frac{d\rho_2}{d\rho_1}(\eta) \frac{d\rho_2}{d\rho_1}(\eta') \rho_1(\eta, \eta')^2 \leq \mu \cdot \lambda \cdot 1 = 1
$$
so equality holds in the inequalities above, hence $\frac{d\rho_2}{d\rho_1}(\eta') = \lambda$ and $\rho_1(\eta, \eta') = 1$. $\diamond$

\medskip

Let $f : (Z_1, \rho_1) \to (Z_2, \rho_2)$ be a homeomorphism between metric spaces. We say $f$ is Moebius if $f$ preserves cross-ratios
with respect to the metrics $\rho_1$ and $\rho_2$, i.e. $[f(\xi), f(\xi'), f(\eta), f(\eta')]_{\rho_2} = [\xi, \xi' ,\eta, \eta']_{\rho_1}$ for all
quadruples of distinct points $\xi,\xi', \eta, \eta'$ in $Z_1$. Then the metrics $\rho_1$ and
$f^* \rho_2$ (the pull-back of $\rho_2$ by $f$) are Moebius equivalent, and we define the derivative of the Moebius map $f$ with respect to
the metrics $\rho_1, \rho_2$ to be the function $\frac{df^* \rho_2}{d\rho_1}$.

\medskip

Let $X$ be a proper, geodesically complete CAT(-1) space (this means that every finite geodesic segment in $X$ can be extended to a
bi-infinite geodesic), with boundary at infinity $\partial X$. The Busemann function of $X$ is the function $B : X \times X \times \partial X \to \mathbb{R}$
defined by
$$
B(x,y,\xi) := \lim_{z \to \xi} (d(x, z) - d(y,z)) \ , \ x,y \in X, \xi \in \partial X
$$
Note that $|B(x,y,\xi)| \leq d(x, y)$ for all $x,y \in X, \xi \in \partial X$. For $x \in X$ and $\xi, \eta \in \partial X$, we denote by $[x, \xi) \subset X$ the unique
geodesic ray joining $x$ to $\xi$, and we denote by $(\xi, \eta) \subset X$ the unique bi-infinite geodesic joining $\xi$ and $\eta$.
For every $x \in X$, there is a metric $\rho_x$ on $\partial X$ called the visual metric on $\partial X$ based at $X$, defined by $\rho_x(\xi, \eta) = e^{-(\xi|\eta)_x}$,
where $(\xi|\eta)_x$ is the Gromov inner product between $\xi,\eta \in \partial X$ with respect to the basepoint $x \in X$, defined by
$$
(\xi|\eta)_x := \lim_{y \to \xi, z \to \eta} \frac{1}{2} (d(x,y)+d(x,z)-d(y,z))
$$
The metric space $(\partial X, \rho_x)$ is compact of diameter one, and the metric $\rho_x$ is antipodal. We have $\rho_x(\xi, \eta) = 1$ if and only if
the point $x$ lies on the bi-infinite geodesic $(\xi, \eta)$.
Moreover, any two visual metrics $\rho_x, \rho_y$ on
$\partial X$ are Moebius equivalent, so there is a canonical cross-ratio function on quadruples of distinct points in $\partial X$, which we will denote by simply $[.,.,.,.]$.
The derivative $\frac{d\rho_y}{d\rho_x}$ is given by
$$
\frac{d\rho_y}{d\rho_x}(\xi) = e^{B(x,y,\xi)}
$$
The space $\mathcal{M}(\partial X, \rho_x)$ is independent of the choice of $x \in X$, and we will denote it by $\mathcal{M}(\partial X)$. The map
$i_X : X \to \mathcal{M}(\partial X), x \mapsto \rho_x$, is an isometric embedding, and the image is $\frac{1}{2}\log 2$-dense in $\mathcal{M}(\partial X)$.

\medskip

For $x \in X$ and a subset $B \subset X$, we define the shadow of the set $B$ as seen from $x$ to be the subset of $\partial X$ defined by
$$
\mathcal{O}(x, B) := \{ \xi \in \partial X | \ [x, \xi) \cap B \neq \emptyset \}
$$

The following lemma will be useful:

\medskip

\begin{lemma} \label{smallshadow} Let $x_0 \in X$ and $R > 0$. For $x \in X$, the diameter of the shadow $\mathcal{O}(x, B(x_0, R))$ with respect to
the visual metric $\rho_x$ tends to $0$ as $x \to \infty$. More precisely, for all $\xi, \eta \in \mathcal{O}(x, B(x_0, R))$,
$$
\rho_x(\xi, \eta) \leq e^{2R -d(x,x_0)}
$$
\end{lemma}

\medskip

\noindent{\bf Proof:} Given $x \in X$ and $\xi \in \mathcal{O}(x, B(x_0, R))$, by definition there exists
$z \in [x, \xi) \cap B(x_0, R)$. Then we have
\begin{align*}
B(x, x_0, \xi) & = B(x, z, \xi) + B(z, x_0, \xi) \\
                & \geq d(x, z) - d(z, x_0) \\
                & \geq d(x, x_0) - 2R \\
\end{align*}
Thus for $\xi, \eta \in \mathcal{O}(x, B(x_0, R))$ we have
\begin{align*}
\rho_x(\xi, \eta)^2 & = \frac{d\rho_x}{d\rho_{x_0}}(\xi)\frac{d\rho_x}{d\rho_{x_0}}(\eta) \rho_{x_0}(\xi, \eta)^2 \\
                    & = e^{B(x_0, x, \xi)} e^{B(x_0, x, \eta)} \rho_{x_0}(\xi, \eta)^2 \\
                    & \leq e^{2R - d(x, x_0)} e^{2R - d(x,x_0)} \\
\end{align*}
and so
$$
\rho_x(\xi, \eta) \leq e^{2R - d(x,x_0)}
$$
$\diamond$

\medskip

The space of geodesics $\mathcal{G}X$ of $X$ is defined to be the space $\mathcal{G}X := \{ \gamma : \mathbb{R} \to X | \ \gamma \hbox{ is an isometric embedding } \}$
equipped with the topology of uniform convergence on compacts. We define continuous maps $\pi : \mathcal{G}X \to X$ and $p : \mathcal{G}X \to \partial X$ by
$\pi(\gamma) = \gamma(0) \in X$ and $p(\gamma) = \gamma(+\infty) \in \partial X$, and
for $x \in X$, we define $T^1_x X := \pi^{-1}(x) \subset \mathcal{G}X$.
The geodesic flow of the CAT(-1) space $X$ is the one-parameter group of homeomorphisms
$(\phi_t : \mathcal{G}X \to \mathcal{G}X)_{t \in \mathbb{R}}$ defined by $(\phi_t \gamma)(s) := \gamma(s+t)$. When $X$ is a simply connected, complete
Riemannian manifold of negative sectional curvature $K \leq -1$, then the map $\mathcal{G}X \to T^1 X, \gamma \mapsto \gamma'(0)$ is a
homeomorphism conjugating the geodesic flow on $\mathcal{G}X$ to the usual geodesic flow on $T^1 X$.

\medskip

Let $Y$ be another proper, geodesically complete CAT(-1) space, and suppose there is a Moebius homeomorphism $f : \partial X \to \partial Y$. The Moebius map
$f$ induces a homeomorphism $\phi : \mathcal{G}X \to \mathcal{G}Y$ conjugating the geodesic flows, which is defined as follows: given $\gamma \in \mathcal{G}X$,
let $x = \gamma(0), \xi = \gamma(+\infty), \eta = \gamma(-\infty)$, then $\phi(\gamma)$ is defined to be the unique $\tilde{\gamma} \in \mathcal{G}Y$ such that
$\tilde{\gamma}(+\infty) = f(\xi), \tilde{\gamma}(-\infty) = f(\eta)$, and $\tilde{\gamma}(0) = y$, where $y$ is the unique point in the bi-infinite geodesic
$(f(\eta), f(\xi)) \subset Y$ such that $\frac{df^*\rho_y}{d\rho_x}(\xi) = 1$.

\medskip

In a CAT(-1) space, any bounded set $B \subset X$ has a unique circumcenter $c(B) \in X$,
i.e. the unique point minimizing the function $x \in X \mapsto \sup_{y \in B} d(x,y)$.
For a compact set $K \subset \mathcal{G}X$ such that $p(K) \subset \partial X$ has at least two points, the limit of the circumcenters
$c(\pi(\phi_t(K)))$ exists as $t \to +\infty$, we call the limit the asymptotic circumcenter of the set $K$ and denote it by $c_{\infty}(K) \in X$.
The geodesic conjugacy $\phi : \mathcal{G}X \to \mathcal{G}Y$ induced by a Moebius map $f : \partial X \to \partial Y$ then allows us to define an extension
$F : X \to Y$ of $f$, called the circumcenter extension of $f$, by
$$
F(x) := c_{\infty}(\phi(T^1_x X)) \in Y
$$
The circumcenter extension is a $(1, \log 2)$-quasi-isometry and is
locally $1/2$-Holder. For $x \in X$, the point $F(x) \in Y$ can be characterized as the unique point in $Y$ minimizing the
function $y \in Y \mapsto d_{\mathcal{M}}(f_* \rho_x, \rho_y)$ (where $f_* \rho_x \in \mathcal{M}(\partial Y)$ is the push-forward of $\rho_x \in \mathcal{M}(\partial X)$
by the Moebius map $f$).

\medskip

\section{Some properties of the circumcenter extension}

\medskip

Throughout this section, $X, Y$ will denote two complete, simply connected manifolds with pinched negative curvature $-b^2 \leq K \leq -1$.
Suppose there is a Moebius homeomorphism $f : \partial X \to \partial Y$ with inverse $g : \partial Y \to \partial X$,
and let $F : X \to Y$ and $G : Y \to X$ be the circumcenter extensions of $f$ and $g$ respectively. Then from \cite{biswas6}, we have that
$F$ and $G$ are $\sqrt{b}$-bi-Lipschitz homeomorphisms which are inverses of each other. Define a function $r : X \to \mathbb{R}$ by
$$
r(x) := d_{\mathcal{M}}( f_* \rho_x, \rho_{F(x)}) = \sup_{\xi \in \partial X} \log \frac{df_* \rho_x}{d\rho_{F(x)}}(f(\xi))
$$

In the following, we identify $\mathcal{G}X, \mathcal{G}Y$ with $T^1 X, T^1 Y$ respectively, and we identify the geodesic
conjugacy $\phi : \mathcal{G}X \to \mathcal{G}Y$ with a geodesic conjugacy $\phi : T^1 X \to T^1 Y$. We also identify the maps
$\pi : \mathcal{G}X \to X, p : \mathcal{G}X \to \partial X$ with maps $\pi : T^1 X \to X, p : T^1 X \to \partial X$ respectively
(and similarly for the corresponding maps for $Y$). For $x \in X, \xi \in \partial X$ we denote by
$\overrightarrow{x\xi} \in T^1_x X$ the unit tangent vector at $x$ given by $\gamma'(0)$, where $\gamma$ is the unique geodesic satisfying $\gamma(0) = x,
\gamma(+\infty) = \xi$. The flip $T^1_x X \to T^1_x X, v \mapsto -v$, induces a continuous involution $i_x : \partial X \to \partial X$, defined by
requiring that $\overrightarrow{xi_x(\xi)} = - \overrightarrow{x\xi}$ for all $\xi \in \partial X$. Similarly for $y \in Y$ we have an involution
$i_y : \partial Y \to \partial Y$. The following lemma follows from Lemma 4.13 of \cite{biswas6}:

\medskip

\begin{lemma} \label{derivformula} For $x \in X, y \in Y, \xi \in \partial X$, we have
$$
\log \frac{df_* \rho_x}{d\rho_y}(f(\xi)) = B(y, \pi(\phi(\overrightarrow{x\xi})), f(\xi))
$$
In particular,
$$
r(x) = \sup_{\xi \in \partial X} B(F(x), \pi(\phi(\overrightarrow{x\xi})), f(\xi))
$$
\end{lemma}

\medskip

\begin{lemma} \label{rlipschitz} The function $r : X \to \mathbb{R}$ is $1$-Lipschitz.
\end{lemma}

\medskip

\noindent{\bf Proof:} Let $x, y \in X$. Since $\phi : T^1 X \to T^1 Y$ conjugates the geodesic flows, we have, for any $\xi \in \partial X$,
$$
B(\pi(\phi(\overrightarrow{x\xi})), \pi(\phi(\overrightarrow{y\xi})), f(\xi)) = B(x, y, \xi)
$$
We then have, using Lemma \ref{derivformula} above,
\begin{align*}
r(x) = d_{\mathcal{M}}(f_* \rho_x, \rho_{F(x)}) & \leq d_{\mathcal{M}}( f_* \rho_x, \rho_{F(y)}) \\
                                                & = \sup_{\xi \in \partial X} \log \frac{df_* \rho_x}{d\rho_{F(y)}}(f(\xi)) \\
                                                & = \sup_{\xi \in \partial X} B(F(y), \pi(\phi(\overrightarrow{x\xi})), f(\xi)) \\
                                                &  = \sup_{\xi \in \partial X} B(F(y), \pi(\phi(\overrightarrow{y\xi})), f(\xi))
                                                + B(\pi(\phi(\overrightarrow{y\xi})), \pi(\phi(\overrightarrow{x\xi})), f(\xi) \\
                                                & = \sup_{\xi \in \partial X} B(F(y), \pi(\phi(\overrightarrow{y\xi})), f(\xi)) + B(y, x, \xi) \\
                                                & \leq \sup_{\xi \in \partial X} B(F(y), \pi(\phi(\overrightarrow{y\xi})), f(\xi)) + d(x,y) \\
                                                & = r(y) + d(x, y) \\
\end{align*}

Thus $r(x) - r(y) \leq d(x,y)$. Interchanging $x$ and $y$ the same argument as above gives $r(y) - r(x) \leq d(x,y)$, hence $|r(x) - r(y)| \leq d(x,y)$. $\diamond$

\medskip

We say that a probability measure $\mu$ on $\partial X$ is {\it balanced} at a point $x \in X$ if the vector-valued integral
$\int_{\partial X} \overrightarrow{x\xi} d\mu(\xi) \in T_x X$ equals $0 \in T_x X$, or equivalently if $\int_{\partial X} < v, \overrightarrow{x\xi} > d\mu(\xi) = 0$
for all $v \in T_x X$. If the compact $K \subset \partial X$ denotes the support of $\mu$, then it is shown in \cite{biswas6} that $\mu$ is balanced at $x$ if and only if
the convex hull in $T_x X$ of the compact set $\{ \overrightarrow{x\xi} : \xi \in K \}$ contains the origin of $T_x X$.

\medskip

For $x \in X$, let $K_x \subset \partial X$ denote the set on which the function $\xi \in \partial X \mapsto \frac{df_* \rho_x}{d\rho_{F(x)}}(f(\xi))$
attains its maximum value. In \cite{biswas6}, it is shown that for any $x \in X$, there exists a probability measure $\mu_x$ on $\partial X$ with support
contained in $K_x$ such that the measure $\mu_x$ is balanced at $x$, and such that the measure $f_* \mu_x$ on $\partial Y$ is balanced at $F(x) \in Y$
(with a similar definition of balanced measures for measures on $\partial Y$ and points of $Y$).

\medskip

The main result of this section is the following:

\medskip

\begin{theorem} \label{rconstant} The function $r$ is constant.
\end{theorem}

\medskip

\noindent{\bf Proof:} Since the function $r$ and the circumcenter map $F$ are both Lipschitz, they are differentiable almost everywhere, so the set of points $D \subset X$
at which both $r$ and $F$ are differentiable has full measure. Let $x \in D$ and let $\xi \in K_x$. Then for any $y \in X$,

\begin{align*}
r(y) & \geq B(F(y), \pi(\phi(\overrightarrow{y\xi})), f(\xi)) \\
     & = B(F(y), F(x), f(\xi)) + B(F(x), \pi(\phi(\overrightarrow{x\xi})), f(\xi)) + B(\pi(\phi(\overrightarrow{x\xi})), \pi(\phi(\overrightarrow{y\xi})), f(\xi)) \\
     & = B(F(y), F(x), f(\xi)) + r(x) + B(x, y, \xi) \\
\end{align*}

so
\begin{equation} \label{randb}
r(y) - r(x) \geq B(F(y), F(x), f(\xi)) + B(x, y, \xi)
\end{equation}

for all $y \in X, \xi \in K_x$. It is well-known that the gradient at $x$ of the function $y \in X \mapsto B(x, y, \xi)$ is given by the vector $\overrightarrow{x\xi}$,
while the gradient at $F(x)$ of the function $z \in Y \mapsto B(z, F(x), f(\xi))$ is given by the vector $-\overrightarrow{F(x)f(\xi)}$.
Let $v \in T_x X$ and $t > 0$, and let $y = \exp_x tv \in X$. Then as $t \to 0$, using the fact that $r$ and $F$ are differentiable at $x$, equation (\ref{randb}) above
gives
$$
dr_x(tv) + o(t) \geq -<dF_x(tv), \overrightarrow{F(x)f(\xi)} > + < tv, \overrightarrow{x\xi} > + o(t)
$$
so dividing by $t$ above and letting $t$ tend to $0$ gives
\begin{equation} \label{drdF}
dr_x(v) \geq < v, \overrightarrow{x\xi} > - < dF_x(v), \overrightarrow{F(x)f(\xi)} >
\end{equation}
for all $v \in T_x X$, $\xi \in K_x$. Integrating both sides of inequality (\ref{drdF}) above over the set $K_x$ with respect to the probability measure $\mu_x$,
and using the facts that the support of $\mu_x$ is contained in $K_x$, the measure $\mu_x$ is balanced at $x$ and the measure $f_* \mu_x$ is
balanced at $F(x)$, we obtain
\begin{align*}
dr_x(v) = \int_{\partial X} dr_x(v) d\mu_x(\xi) & = \int_{K_x} dr_x(v) d\mu_x(\xi) \\
                                                & \geq \int_{K_x} < v, \overrightarrow{x\xi} > d\mu_x(\xi) - \int_{K_x} < dF_x(v), \overrightarrow{F(x)f(\xi)} > d\mu_x(\xi) \\
                                                & = \int_{\partial X} < v, \overrightarrow{x\xi} > d\mu_x(\xi) - \int_{\partial X} < dF_x(v), \overrightarrow{F(x)f(\xi)} > d\mu_x(\xi) \\
                                                & = \int_{\partial X} < v, \overrightarrow{x\xi} > d\mu_x(\xi) - \int_{\partial Y} < dF_x(v), \overrightarrow{F(x)\eta} > df_* \mu_x(\eta) \\
                                                & = 0 \\
\end{align*}

Thus $dr_x(v) \geq 0$ for all $v \in T_x X$, replacing $v$ by $-v$ gives $dr_x(-v) \geq 0$ so $dr_x(v) \leq 0$ for all $v \in T_x X$, and hence $dr_x(v) = 0$ for all $v \in T_x X$.
Since $r$ is Lipschitz and $dr_x = 0$ for $x$ in the full measure set $D$, it follows that $r$ is constant. $\diamond$

\medskip

A corollary of the proof of the above theorem is the following:

\medskip

\begin{prop} \label{dFKx} Let $x \in X$ be a point of differentiability of $F$. Then for all $\xi \in K_x, v \in T_x X$ we have
$$
< dF_x(v), \overrightarrow{F(x)f(\xi)} > = < v, \overrightarrow{x\xi} >
$$
Equivalently,
$$
dF^*_x(\overrightarrow{F(x)f(\xi)}) = \overrightarrow{x\xi}
$$
for all $\xi \in K_x$.
\end{prop}

\medskip

\noindent{\bf Proof:} By the previous theorem the function $r$ is constant, so the set $D$ in the proof of the previous theorem is just the
set of points of differentiability of $F$. Let $x \in D$, and $\xi \in K_x$. Since $r$ is constant, equation (\ref{drdF}) above
gives
$$
0 \geq < v, \overrightarrow{x\xi} > - < dF_x(v), \overrightarrow{F(x)f(\xi)} >
$$
for all $v \in T_x X$. Replacing $v$ by $-v$ in the above equation gives
$$
0 \leq < v, \overrightarrow{x\xi} > - < dF_x(v), \overrightarrow{F(x)f(\xi)} >
$$
for all $v \in T_x X$. Combining the two gives $< dF_x(v), \overrightarrow{F(x)f(\xi)} > = < v, \overrightarrow{x\xi} >$ for all $v \in T_x X$. $\diamond$

\medskip

\begin{lemma} \label{qisom} Let $M \geq 0$ denote the constant value of the function $r$. Then the circumcenter map $F : X \to Y$ is a $(1, 2M)$-quasi-isometry, i.e.
$$
d(x,y) - 2M \leq d(F(x), F(y)) \leq d(x,y) + 2M
$$
for all $x,y \in X$.
\end{lemma}

\medskip

\noindent{\bf Proof:} Note that push-forward of metrics by $f$ gives an isometry $f_* : \mathcal{M}(\partial X) \to \mathcal{M}(\partial Y)$. So for $x,y \in X$,
we have
\begin{align*}
d(x,y) & = d_{\mathcal{M}}( \rho_x, \rho_y) \\
       & = d_{\mathcal{M}}( f_* \rho_x, f_* \rho_y) \\
       & \leq d_{\mathcal{M}}( f_* \rho_x, \rho_{F(x)}) + d_{\mathcal{M}}( \rho_{F(x)}, \rho_{F(y)}) + d_{\mathcal{M}}( \rho_{F(y)}, f_* \rho_y) \\
       & = M + d(F(x), F(y)) + M \\
\end{align*}
Similarly,
\begin{align*}
d(F(x), F(y)) & = d_{\mathcal{M}}( \rho_{F(x)}, \rho_{F(y)}) \\
       & \leq d_{\mathcal{M}}( \rho_{F(x)}, f_* \rho_x) + d_{\mathcal{M}}( f_* \rho_x, f_* \rho_y) + d_{\mathcal{M}}( f_* \rho_y, \rho_{F(y)}) \\
       & = M + d(x, y) + M \\
\end{align*}
thus
$$
d(x,y) - 2M \leq d(F(x), F(y)) \leq d(x,y) + 2M
$$
$\diamond$

\medskip

The following lemma is a straightforward consequence of Lemma \ref{maxminantipodal}:

\medskip

\begin{lemma} \label{maxminflip} Let $x \in X$ and $y \in Y$. Then:

\medskip

\noindent (i) The function $\frac{d\rho_x}{df^*\rho_y}$ attains its maximum at $\xi \in \partial X$ if and only if it attains its minimum at $i_x(\xi)$.
Moreover in this case $f(i_x(\xi)) = i_y(f(\xi))$, so $y$ lies on the bi-infinite geodesic $(f(\xi), f(i_x(\xi)))$.

\medskip

\noindent (ii) If $\xi \in \partial X$ is a maximum of $\frac{d\rho_x}{df^*\rho_y}$ then the point $z = \pi(\phi(\overrightarrow{x\xi})) \in Y$ is the unique
point on the geodesic ray $[y, f(\xi)) \subset Y$ at a distance $d_{\mathcal{M}}(f_* \rho_x, \rho_y)$ from $y$.

\end{lemma}

\medskip

\noindent{\bf Proof:} (i) We first note that since $X$ is a simply connected manifold of
negative curvature, for $\xi, \eta \in \partial X$ we have $\rho_x(\xi, \eta) = 1$ if and only if $\eta = i_x(\xi)$.
Let $\xi \in \partial X$ be a maximum of $\frac{d\rho_x}{df^*\rho_y}$. Let $\eta = f^{-1}(i_y(f(\xi))) \in \partial X$, then
$f^* \rho_y(\xi, \eta) = \rho_y(f(\xi), f(\eta)) = \rho_y(f(\xi), i_y(f(\xi))) = 1$, hence by Lemma \ref{maxminantipodal} we have that
$\eta$ is a minimum of $\frac{d\rho_x}{df^*\rho_y}$. Moreover, by Lemma \ref{maxminantipodal}, $\rho_x(\xi, \eta) = 1$, thus $\eta = i_x(\xi)$, so
$\frac{d\rho_x}{df^*\rho_y}$ attains its minimum at $i_x(\xi)$, and $f(i_x(\xi)) = i_y(f(\xi))$.

\medskip

For the converse, suppose that $i_x(\xi) \in \partial X$ is a minimum of $\frac{d\rho_x}{df^*\rho_y}$. Then $\rho_x(\xi, i_x(\xi)) = 1$ implies by Lemma
\ref{maxminantipodal} that $\frac{d\rho_x}{df^*\rho_y}$ attains its maximum at $\xi$. Moreover, by Lemma \ref{maxminantipodal},
$f^* \rho_y(\xi, i_x(\xi)) = 1$, so $\rho_y(f(\xi), f(i_x(\xi))) = 1$, hence $f(i_x(\xi)) = i_y(f(\xi))$.

\medskip

\noindent (ii) Let $\xi$ be a maximum of $\frac{d\rho_x}{df^*\rho_y}$. By definition of the geodesic conjugacy $\phi$, the point $z = \pi(\phi(\overrightarrow{x\xi})) \in Y$
lies on the bi-infinite geodesic $(f(\xi), f(i_x(\xi))) \subset Y$. By (i) above, the point $y$ also lies on the bi-infinite geodesic $(f(\xi), f(i_x(\xi)))$.
Since $\xi$ is a maximum of $\frac{d\rho_x}{df^*\rho_y}$ it follows that
$\log \frac{d\rho_x}{df^*\rho_y}(\xi) = d_{\mathcal{M}}( \rho_x, f^*\rho_y) = d_{\mathcal{M}}( f_* \rho_x, \rho_y)$ (note that push-forward of metrics by $f$ gives
an isometry $f_* : \mathcal{M}(\partial X) \to \mathcal{M}(\partial Y)$). Thus by Lemma \ref{derivformula} we have
\begin{align*}
B(y, z, f(\xi)) & = \log \frac{df_*\rho_x}{d\rho_y}(f(\xi)) \\
                & = \log \frac{d\rho_x}{df^*\rho_y}(\xi) \\
                & = d_{\mathcal{M}}( f_* \rho_x, \rho_y) \\
\end{align*}
Since $y,z$ both lie on the geodesic $(f(\xi), f(i_x(\xi)))$, it follows that $z$ is the unique point on the geodesic ray $[y, f(\xi))$ at a distance
$d_{\mathcal{M}}( f_* \rho_x, \rho_y)$ from $y$. $\diamond$

\medskip

Finally, we need a lemma about Riemannian angles and comparison angles from \cite{biswas5}.
For $x \in X$ and $\xi, \eta \in \partial X$, let $\angle \xi x \eta \in [0, \pi]$
denote the Riemannian angle between the geodesic rays $[x, \xi)$ and $[x, \eta)$ at the point $x$. Then the following holds (this is Lemma 6.6
of \cite{biswas5}):

\medskip

\begin{lemma} \label{anglecomp} For all $x \in X$ and $\xi, \eta \in \partial X$ we have
$$
\rho_x(\xi, \eta)^b \leq \sin\left(\frac{1}{2}\angle \xi x \eta \right) \leq \rho_x(\xi, \eta)
$$
\end{lemma}

\medskip

\section{Proof of main theorem}

\medskip

Let $(X, g_0)$ be a complete, simply connected manifold of pinched negative curvature $-b^2 \leq K_{g_0} \leq -1$. Let $g_1$ be a metric on $X$
such that $g_1 = g_0$ outside a compact in $X$, and suppose $g_1$ is negatively curved, $K_{g_1} \leq -1$. Then the metrics $g_0, g_1$ are bi-Lipschitz,
so the identity map $id : (X, g_0) \to (X, g_1)$ induces a homeomorphism between boundaries which we denote by
$f : \partial_{g_0} X \to \partial_{g_1} X$. Suppose the map $f$ is Moebius. Let $F : (X, g_0) \to (X, g_1)$ be the circumcenter extension of the Moebius map $f$.
Note that both metrics $g_0, g_1$ have pinched negative curvature (since $g_0$ does, and $g_1 = g_0$ outside a compact), so the results of the previous
section apply to $F$. In particular, by Theorem \ref{rconstant}, the function $r(x) = d_{\mathcal{M}}( f_* \rho_x, \rho_{F(x)} )$ is constant, let
$M \geq 0$ denote its constant value. By Lemma \ref{qisom}, to show that the circumcenter map $F$ is an isometry, it suffices to show that $M = 0$.

\medskip

Let $T^1 X_{g_0} \subset TX$ and $T^1 X_{g_1} \subset TX$ denote the unit tangent bundles with respect to the metrics $g_0, g_1$ respectively, and let
$\phi : T^1 X_{g_0} \to T^1 X_{g_1}$ denote the geodesic conjugacy induced by the Moebius map $f$. For $x \in X$, let $\rho^{g_0}_x$ and $\rho^{g_1}_x$
denote the visual metrics based at $x$ on the boundaries $\partial_{g_0} X$ and $\partial_{g_1} X$ of $(X, g_0)$ and $(X, g_1)$ respectively.
For $x \in X$ and $\xi, \eta \in \partial_{g_i} X$, let $(\xi, \eta)_i \subset X$ denote the bi-infinite $g_i$-geodesic with endpoints $\xi, \eta$, and
let $[x, \xi)_i \subset X$ denote the $g_i$-geodesic ray joining $x$ to $\xi$, and let $\overrightarrow{x\xi}^i \in T^1_x X_{g_i}$
denote the $g_i$-unit tangent vector to the $g_i$-geodesic ray $[x, \xi)_i$ at the point $x$, where $i = 0,1$.
For $x \in X$ and a compact $K \subset X$, let $\mathcal{O}_i(x, K) \subset \partial_{g_i} X$ denote the shadow of the set $K$ as seen from the point $x$
with respect to the metric $g_i$, where $i = 0,1$.
For $i = 0,1$ and $x \in X$, let $i^{g_i}_x : \partial_{g_i} X \to \partial_{g_i} X$ denote the involution of the boundary of $(X, g_i)$
as defined in the previous section.

\medskip

\begin{lemma} \label{conjid} Let $K = $ supp$(g_1 - g_0)$ denote the support of the symmetric 2-tensor $g_1 - g_0$. Let $x \in X - K$.
If $\xi \in \partial_{g_0} X$ is such that $\xi, i^{g_0}_x(\xi) \in \partial_{g_0} X - \mathcal{O}_o(x, K)$,
then $\overrightarrow{x\xi}^0 = \overrightarrow{xf(\xi)}^1 \in T^1_x X_{g_0} \cap T^1_x X_{g_1}$ and
$\phi(\overrightarrow{x\xi}^0) = \overrightarrow{xf(\xi)}^1 = \overrightarrow{x\xi}^0$.
\end{lemma}

\medskip

\noindent{\bf Proof:} The hypothesis on $\xi$ implies that the $g_0$-geodesic rays $[x, \xi)_0$ and $[x, i^{g_0}_x(\xi))_0$ are disjoint from $K$,
hence so is the bi-infinite $g_0$-geodesic $(\xi, i^{g_0}_x(\xi))_0$, thus it is also a $g_1$-geodesic, hence $(\xi, i^{g_0}_x(\xi))_0$ equals the
bi-infinite $g_1$-geodesic $(f(\xi), f(i^{g_0}_x(\xi)))_1$.
In particular $\overrightarrow{x\xi}^0 = \overrightarrow{xf(\xi)}^1 \in T^1_x X_{g_0} \cap T^1_x X_{g_1}$,
and $\phi(\overrightarrow{x\xi}^0)$ is tangent to $(\xi, i^{g_0}_x(\xi))_0$, so $\pi(\phi(\overrightarrow{x\xi}^0))$ lies on $(\xi, i^{g_0}_x(\xi))_0$.
Now we can choose a neighbourhood $U$ of $\xi$ in $\partial_{g_0} X$ which is disjoint from $\mathcal{O}_o(x, K)$, and such that for any $\eta \in U$, the $g_0$-geodesic
$(\xi, \eta)_0$ is disjoint from $K$ (by choosing $U$ small enough). Then for $\eta \in U$, the $g_0$-geodesics $[x, \xi)_0, [x, \eta)_0, (\xi, \eta)_0$ are
disjoint from $K$, hence they are $g_1$-geodesics as well, and it follows that $\rho^{g_0}_x(\xi, \eta) = \rho^{g_1}_x(f(\xi), f(\eta))$ for all $\eta \in U$.
Hence
$$
\frac{df^*\rho^{g_1}_x}{\rho^{g_0}_x}(\xi) = \lim_{\eta \to \xi} \frac{f^*\rho^{g_1}_x(\xi, \eta)}{\rho^{g_0}_x(\xi, \eta)} = 1
$$
so it follows from the definition of $\phi$ that $\pi(\phi(\overrightarrow{x\xi}^0)) = x$, thus
$\phi(\overrightarrow{x\xi}^0) = \overrightarrow{xf(\xi)}^1 = \overrightarrow{x\xi}^0$. $\diamond$

\medskip

For $i = 0,1$, let $d_{g_i}$ denote the distance function of $(X, g_i)$, and for $x \in X$ and $\xi, \eta \in \partial_{g_i} X$ let $\angle_i \xi x \eta$ denote the Riemannian
angle between the $g_i$-geodesic rays $[x, \xi)_i, [x, \eta)_i$ at the point $x$ with respect to the metric $g_i$.

\medskip

We can now prove the main theorem:

\medskip

\noindent{\bf Proof of Theorem \ref{mainthm}:} As remarked earlier, it suffices to show that the constant $M = 0$,
where $d_{\mathcal{M}}( f_* \rho^{g_0}_x, \rho^{g_1}_{F(x)} ) = M$ for all $x \in X$. Fix $\epsilon > 0$, we will show that $M \leq \epsilon$.

\medskip

Fix a basepoint $x_0 \in X$ and choose $R > 0$
such that the support of $g_1 - g_0$ is contained in the $g_0$-ball of radius $R$ around $x_0$, and
let $B$ denote the closed $g_0$-ball of radius $R$ around $x_0$. Fix $\xi_0, \eta_0 \in \partial_{g_0} X$ such that $x_0 \in (\xi_0, \eta_0)_0$,
let $\gamma : \mathbb{R} \to X$ be the unique unit speed $g_0$-geodesic such that $\gamma(-\infty) = \xi_0, \gamma(0) = x_0, \gamma(+\infty) = \eta_0$.
For $t > R$ let $x_t \in X$ denote the point $\gamma(t)$, and define $\epsilon_t > 0$ by
$$
\epsilon_t := \sup \{ \angle_0 \xi x_t \xi_0 | \xi \in \mathcal{O}_0(x_t, B) \}
$$

Then it follows from Lemma \ref{smallshadow} and Lemma \ref{anglecomp} that $\epsilon_t \to 0$ as $t \to +\infty$.

\medskip

Let $K_t \subset \partial_{g_0} X$ denote the set where the function $\frac{d\rho^{g_0}_{x_t}}{df^*\rho^{g_1}_{F(x_t)}}$ attains its maximum value $e^M$. 
Let $C_t \subset T_{x_t} X$ denote the cone
$$
C_t := \{ v \in T_{x_t} X | < v, \overrightarrow{x_t \xi_0}^0 >_{g_0} \geq \cos(\epsilon_t) ||v||_{g_0} \}
$$
and let $D_t := \{ -v \in T_{x_t} X | v \in C_t \}$. Then for $\xi \in \partial_{g_0} X$, if $\overrightarrow{x_t \xi}^0 \notin C_t \cup D_t$,
 then $\xi, i^{g_0}_{x_t}(\xi) \notin \mathcal{O}_0(x_t, B)$.
 Moreover, for $v, w \in C_t$
and $\alpha, \beta \geq 0$ we have $\alpha v + \beta w \in C_t$. Now if $\xi, \eta \in \partial_{g_0} X$ are such that $\overrightarrow{x\xi}^0 \in C_t$ 
and $\overrightarrow{x\eta}^0 \in D_t$, then by the triangle inequality
$$
\rho^{g_0}_{x_t}(\xi, \eta) \geq 1 - \rho^{g_0}_{x_t}(\xi, \xi_0) - \rho^{g_0}_{x_t}(\eta, \eta_0)
$$
and by Lemma \ref{anglecomp} we have
$$
\rho^{g_0}_{x_t}(\xi, \xi_0)^b \leq \sin(\epsilon_t/2), \rho^{g_0}_{x_t}(\eta, \eta_0)^b \leq \sin(\epsilon_t/2),
$$
so since $\epsilon_t \to 0$ as $t \to +\infty$, by choosing $t > R$ large enough we may assume that
$$
\rho^{g_0}_{x_t}(\xi, \eta) \geq e^{-\epsilon}
$$
whenever $\xi, \eta \in \partial_{g_0} X$ are such that $\overrightarrow{x\xi}^0 \in C_t$
and $\overrightarrow{x\eta}^0 \in D_t$. We fix such a $t > R$ large enough so that this holds. 

\medskip

As stated in section 3, there exists a probability measure $\mu$ on $\partial_{g_0} X$ with support contained in $K_t$ such that $\mu$ is balanced
at $x_t \in (X, g_0)$, equivalently the convex hull in $T_{x_t} X$ of the compact set $\{ \overrightarrow{x_t \xi}^0 | \xi \in K_t \}$ contains the
origin of $T_{x_t} X$. By the classical Caratheodory theorem on convex hulls, it follows that there exist distinct points $\xi_1, \dots, \xi_k \in K_t$ and
$\alpha_1, \dots, \alpha_k > 0$ such that $\alpha_1 \overrightarrow{x_t\xi_1}^0 + \dots + \alpha_k \overrightarrow{x_t\xi_k}^0 = 0$ and
$\alpha_1 + \dots + \alpha_k = 1$, where $1 \leq k \leq n+1$ (here $n$ is the dimension of $X$). Note that since the vectors $\overrightarrow{x_t \xi_i}^0$
are nonzero, we must have $k \geq 2$. We now consider various cases:

\medskip

\noindent{\bf Case 1.} $k = 2$: 

\medskip

Then since $\overrightarrow{x_t \xi_1}^0, \overrightarrow{x_t\xi_2}^0$ are unit vectors, the relation
$\alpha_1 \overrightarrow{x_t\xi_1}^0 + \alpha_2 \overrightarrow{x_t\xi_2}^0 = 0$ implies that
$\overrightarrow{x_t \xi_1}^0 = - \overrightarrow{x_t \xi_2}^0$, hence $\xi_2 = i^{g_0}_{x_t}(\xi_1)$. By Lemma \ref{maxminflip}, the function
$\frac{d\rho^{g_0}_{x_t}}{df^*\rho^{g_1}_{F(x_t)}}$ attains its minimum at $\xi_2$, so since $\xi_2 \in K_t$, the maximum and minimum of the function
$\frac{d\rho^{g_0}_{x_t}}{df^*\rho^{g_1}_{F(x_t)}}$ are equal, thus $e^M = e^{-M}$, and so $M = 0$ as required.

\medskip

\noindent{\bf Case 2.} $k \geq 3$, and there exist $1 \leq i \neq j \leq k$ such that $\overrightarrow{x_t \xi_i}^0, \overrightarrow{x_t \xi_j}^0 \in T_{x_t} X - (C_t \cup D_t)$: 

\medskip

In this case, $\xi_i, i^{g_0}_{x_t}(\xi_i), \xi_j, i^{g_0}_{x_t}(\xi_j) \in \partial_{g_0} X - \mathcal{O}_o(x_t, B)$.
It follows from Lemma \ref{conjid} that the points
$z_i := \pi(\phi(\overrightarrow{x_t\xi_i}^0)), z_j := \pi(\phi(\overrightarrow{x_t\xi_j}^0))$ satisfy $z_i = x_t = z_j$. Thus the $g_1$-geodesics $(f(\xi_i), f(i^{g_0}_{x_t}(\xi_i)))_1$ and
$(f(\xi_j), f(i^{g_0}_{x_t}(\xi_j)))_1$ intersect at the point $x_t$. On the other hand, by Lemma \ref{maxminflip}, the geodesics
$(f(\xi_i), f(i^{g_0}_{x_t}(\xi_i)))_1$ and
$(f(\xi_j), f(i^{g_0}_{x_t}(\xi_j)))_1$ intersect at the point $F(x_t)$. If $\xi_j \neq i^{g_0}_{x_t}(\xi_i)$, then the geodesics
$(f(\xi_i), f(i^{g_0}_{x_t}(\xi_i)))_1$ and $(f(\xi_j), f(i^{g_0}_{x_t}(\xi_j)))_1$ have a unique point of intersection, thus $x_t = F(x_t)$, and by
Lemma \ref{maxminflip} we have
$$
M = d_{\mathcal{M}}( f_* \rho^{g_0}_{x_t}, \rho^{g_1}_{F(x_t)} ) = d_{g_1}(z_i, F(x_t)) = d_{g_1}(x_t, x_t) = 0
$$
If on the other hand $\xi_j = i^{g_0}_{x_t}(\xi_i)$, then the same argument as in Case 1 above shows that $M = 0$. Thus in either case $M = 0$.

\medskip

\noindent{\bf Case 3.} $k \geq 3$, and $\overrightarrow{x\xi}^0 \in T_{x_t} - (C_t \cup D_t)$ for at most one $i \in \{ 1, \dots, k \}$:

\medskip

Then relabelling the $\xi_i$'s if necessary, we may assume that $\overrightarrow{x_t\xi_1}^0, \dots, \overrightarrow{x_t\xi_{k - 1}}^0 \in C_t \cup D_t$.
Now if $\overrightarrow{x_t\xi_1}^0, \dots, \overrightarrow{x_t\xi_{k - 1}}^0 \in C_t$, then
$\alpha_1 \overrightarrow{x_t\xi_1}^0 + \dots + \alpha_{k-1} \overrightarrow{x_t\xi_{k - 1}}^0 \in C_t$ and it follows that $\overrightarrow{x_t \xi_k}^0 \in D_t$.
 Similarly if $\overrightarrow{x_t\xi_1}^0, \dots, \overrightarrow{x_t\xi_{k - 1}}^0 \in D_t$, then we must have $\overrightarrow{x_t \xi_k}^0 \in C_t$. Thus either way,
 there exist $1 \leq i \neq j \leq k$ such that $\overrightarrow{x_t \xi_i}^0 \in C_t$ and $\overrightarrow{x_t \xi_j}^0 \in D_t$. Let $\eta = i^{g_0}_{x_t}(\xi_i),
 \eta' = i^{g_0}_{x_t}(\xi_j)$, then $\overrightarrow{x_t \eta}^0 \in D_t$ and $\overrightarrow{x_t \eta'}^0 \in C_t$, and by Lemma \ref{maxminflip}, the function
$\frac{d\rho^{g_0}_{x_t}}{df^*\rho^{g_1}_{F(x_t)}}$ attains its minimum value $e^{-M}$ at the points $\eta, \eta'$. Now by our hypothesis on $t$ we have 
$$
\rho^{g_0}_{x_t}( \eta, \eta') \geq e^{-\epsilon}.
$$
We then have
\begin{align*}
e^{-2\epsilon} & \leq \rho^{g_0}_{x_t}( \eta, \eta')^2 \\
               & = \frac{d\rho^{g_0}_{x_t}}{df^*\rho^{g_1}_{F(x_t)}}(\eta) \frac{d\rho^{g_0}_{x_t}}{df^*\rho^{g_1}_{F(x_t)}}(\eta') f^*\rho^{g_1}_{F(x_t)}(\xi, \xi')^2 \\
               & \leq e^{-M} \cdot e^{-M} \cdot 1 \\
\end{align*}
thus $e^{-2\epsilon} \leq e^{-2M}$, hence $M \leq \epsilon$.

\medskip

Since Cases 1,2,3 above exhaust all possibilities, it follows that $M \leq \epsilon$ for any given $\epsilon > 0$, thus $M = 0$ as required. $\diamond$

\bibliography{moeb}
\bibliographystyle{alpha}

\end{document}